\documentclass[11pt,a4paper]{article}
\usepackage{amsmath,amsthm,amssymb}

\newtheorem*{theorem}{Theorem}
\newtheorem*{corollary}{Corollary}
\newtheorem*{proposition}{Proposition}

\newcommand*\mapp[3]{#1\colon#2\to#3}
\newcommand\fbaseset{\mathbb B}
\newcommand\fbase{(\fbaseset,\Omega)}
\newcommand\fpinbase{b\in\fbaseset}
\newcommand*\cl[1]{\mathop{\mathrm{cl}}#1}
\newcommand*\tube[2]{{#1}^{{\gets}}#2}

\newcommand*\opentop[1]{\Omega(#1)}
\newcommand\fbasetop{\Omega}
\newcommand*\Neighbs[2]{\mathcal N_{#2}(#1)}

\newcommand*\Cat[1]{\mathbf{#1}}
\newcommand*\Fun[1]{\mathtt{#1}}

\newcommand*\Supp[2]{\left|#1\right|_{#2}}

\DeclareMathOperator\Img{Im}

\newcommand\arrowsrightleft[3]{\mathop{\vcenter{\baselineskip=2ex
     \hbox to#1
    {\rightarrowfill}
     \kern-1.1ex\hbox to #1{\leftarrowfill}}}\limits^{\ #2\ }_{\ #3\ }}

\begin{document}
\title{On Faded Cosheaves of Sets}
\author{Alexei Zouboff}
\date{October 20, 1998}
\maketitle

\begin{abstract}
We prove that the category of faded cosheaves in $\Cat{Set}$
over a sober topological space $\fbase$ is equivalent to a category
$\mathop{\Cat{Sett}}\fbase$ having
the same class of objects as $\Cat{Set}/\fbaseset$ has, but
generally a wider class of morphisms.
We also prove the converse: if each cosheaf over $T_0$\mbox{-}space
$\fbase$ is isomorphic
to the cosheaf of tubes of an appropriate map ranging in $\fbase$, then
$\fbase$ is sober.
\end{abstract}

Throughout the paper $\fbase$ is a (fixed) topological space having
$\fbaseset$ as its underlying set and $\fbasetop$ its open topology. We
denote $\opentop V=\{\,W\in\fbasetop:W\subset V\,\}$ for each
$V\in\fbasetop$ and $\Neighbs xV=\{\,W\in\opentop V:x\in W\,\}$
for each $V\in\fbasetop$ and $x\in V$.

Given two sets $X,Y$ and maps $\mapp fX\fbaseset$, $\mapp gY\fbaseset$,
let us say that a map $\mapp\phi XY$ is a \emph{tubewise morphism\/}
$\mapp\phi fg$, if $\phi(\tube fV)\subset\tube gV$ whenever
$V\in\fbasetop$ (or, that is the same, any of the following equivalent
conditions holds: (i)~$\tube fV\subset\tube{(g\phi)}V$ for each $V\in
\fbasetop$; (ii)~$\tube{(g\phi)}F\subset\tube fF$ for each closed
$F\subset\fbaseset$; (iii)~$\tube{(g\phi)}{(\cl\{b\})}\subset
\tube f{(\cl\{b\})}$ for each $\fpinbase$). Clearly, any
\emph{fibrewise\/} morphism $\mapp\phi fg$ (that is a map $\phi$ such
that $f=g\phi$) is tubewise, and converse is true, if $\fbase$ is a
$T_1$\nobreakdash-space.

The class of all the maps of sets into $\fbaseset$ and the class of
their tubewise morphisms (with the arbitrary law of composition)
form a category, which is denoted by $\mathop{\Cat{Sett}}\fbase$.
By the remark above, $\mathop{\Cat{Sett}}\fbase$ is a supercategory
of the comma category $\Cat{Set}/\fbaseset$ (having the same class of
objects) and coincides with the latter,
if $\fbase$ is a $T_1$\nobreakdash-space.

Given a map $\mapp fX\fbaseset$, there is a cosheaf $\Fun{CS}(f)$
in $\Cat{Set}$ over
the space $\fbase$, defined by $V\mapsto\tube fV$ ($V\in\fbasetop$) and
obvious inclusion maps. This cosheaf is \emph{faded\/}---i.e., all its
inclusion morphisms are monic ($=$injective). Moreover, the correspondence
$f\mapsto\Fun{CS}(f)$ is extendable to a functor
from $\mathop{\Cat{Sett}}\fbase$ to the category of faded cosheaves
in $\Cat{Set}$ over $\fbase$ (the latter category is denoted by
$\Cat{Coshv_m}\fbase$ below), that sends each tubewise morphism
$\mapp\phi fg$ to the family $\varPhi={\{\phi_V\}}_{V\in\fbasetop}$,
where $\phi_V$ is the restriction of $\phi$ onto domain
$\tube fV$ and codomain $\tube gV$ for each $V\in\fbasetop$
(so that $\mapp\varPhi
{\Fun{CS}(f)}{\Fun{CS}(g)}$ is a morphism of faded cosheaves).

Now suppose that $\fbase$ is a sober space. (Recall that
$\fbase$ is \emph{sober}, or \emph{primal} (see, e.g.,
\cite{JohnSTONE}), if it is naturally homeomorphic to the space of
principal prime ideals of the lattice $\fbasetop$, or, equivalently,
if each \emph{completely prime filter} in the lattice $\fbasetop$
(i.e., inaccessible by arbitrary joins, not only finite ones) coincides
with $\Neighbs b\fbaseset$ for a unique $\fpinbase$.
One has $\{\,T_0$\nobreakdash-spaces$\,\}\supset\{\,$sober spaces$\,\}
\supset\{\,T_2$\nobreakdash-spaces$\,\}$, and sobriety is incomparable
with $T_1$ axiom.)

Let $\boldsymbol X=\{X_V, R^W_V\}$ be a faded cosheaf of sets over
$\fbase$ (where $V\mapsto X_V$ for each $V\in\fbasetop$ and
$\mapp{R^W_V}{X_W}{X_V}$ are inclusion maps for $W\in\opentop V$).
Then we have
\begin{equation}\label{eq:1}
\begin{gathered}
\bigcup\{\,\Img R^W_V:W\in\mathcal W\,\}=
\Img R^{\bigcup\mathcal W}_V
\quad\text{for } V\in\fbasetop\text{ and }
\mathcal W\subset\opentop V,                     \\
\Img R^U_V\cap\Img R^W_V=\Img R^{U\cap W}_V\quad
\text{for }V\in\fbasetop \text{ and } W,U\in\opentop V
\end{gathered}
\end{equation}
(these relations hold for arbitrary
cosheaves of sets, if $\bigcup\mathcal W=V$ and $U\cup W=V$,
respectively).

It follows from~\eqref{eq:1} that for any
$V\in\fbasetop$ and each $x\in X_V$, the family
$\mathcal A_V(x)=\{\,W\in\opentop V:x\in\Img R^W_V\,\}$ is a completely
prime filter in $\opentop V$; therefore we have $\mathcal A_V(x)=
\Neighbs bV$ for a unique $b\in V$ (since sobriety is inherited
by open subsets). We have constructed a map $\mapp{\Supp\cdot V}
{X_V}V$ sending each $x\in X_V$ to the $b\in V$ described before.
So, for any $W\in\opentop V$ we have
\begin{equation}\label{eq:2}
\Supp xV\in W\iff x\in\Img R^W_V
\end{equation}
(in other words, the singleton $\{\Supp xV\}$ is exactly
the \emph{support\/} of $x$---cf.~\cite[ch.~5]{Bredon}). Furthemore,
for $W\in\opentop V$ and $x\in X_W$ we have
$\mathcal A_W(x)=\mathcal A_V(R^W_V(x))\cap\opentop W$, since
$R^W_V$ is injective. It follows that
\begin{equation}\label{eq:3}
\Supp\cdot V\circ R^W_V=j^W_V\circ\Supp\cdot W\quad\text{for
}W\in\opentop V,
\end{equation}
where $\mapp{j^W_V}WV$ is the arbitrary embedding map.
We have almost proved the following
\begin{theorem}
If $\fbase$ is sober, then there exists
an equivalence
$\mathop{\Cat{Sett}}\fbase\allowbreak
\arrowsrightleft{1.8em}{\Fun{CS}}{\Fun L}
\Cat{Coshv_m}\fbase$, such that $\Fun L$ is a left inverse to the
functor $\Fun{CS}$.
\end{theorem}

\begin{proof}
Let $\boldsymbol X=\{\,X_V, R^W_V\,\}$ be a faded cosheaf.
Put $\Fun L(\boldsymbol X)=\Supp\cdot\fbaseset$, where
$\mapp{\Supp\cdot\fbaseset}{X_\fbaseset}\fbaseset$ is the map defined
above by~\eqref{eq:2} for $V=\fbaseset$.

For each morphism
${\varPhi=\allowbreak\{\,\mapp{\phi_V}{X_V}{Y_V}\,\}}_{V\in\fbasetop}$ of
faded cosheaves
$\boldsymbol X=\{\,X_V, R^W_V\,\}$ and $\boldsymbol Y=\{\,Y_V, T^W_V\,\}$
put $\Fun L(\varPhi)=\phi_\fbaseset$.
Using~\eqref{eq:2} and~\eqref{eq:3}, one can easily check that
the functor $\Fun L$ is well-defined and satisfies the
conditions required.
\end{proof}

\begin{corollary}
If $\fbase$ is a sober $T_1$-space, then there exists
an equivalence
$\Cat{Set}/\fbaseset\allowbreak
\arrowsrightleft{1.8em}{\Fun{CS}}{\Fun L}
\Cat{Coshv_m}\fbase$.
\end{corollary}

For the conclusion, we will show that sobriety of the base space
$\fbase$ is essential in our considerations: the Proposition stated below
is a converse (in a certain sence) to the Theorem.

It is easy to see that the copresheaf $\{\,X_V, R^W_V\,\}$
over $\fbase$ defined by
\begin{equation}\label{eq:4}
V\mapsto X_V=\{\,\text{all the completely prime filters
           in $\opentop V$}\,\},
\end{equation}
with the injective inclusion maps $\mapp{R^W_V}{X_W}{X_V}$ defined for
$W\subset V$ by
\begin{equation}\label{eq:5}
\mathcal A\overset{R^W_V}\longmapsto
       \{\,U\in\opentop V:U\cap W\in\mathcal A\,\},
\end{equation}
is actually a faded cosheaf of sets.
(Proof. Let $\mathcal W\subset \fbasetop$, $V=\bigcup\mathcal W$ and
${\{\,\mapp{\phi_W}{X_W}M\,\}}_{W\in\mathcal W}$ a sink of maps such that
$\phi_U\circ R^{U\cap W}_U=\phi_W\circ R^{U\cap W}_W$ whenever
$U,W\in\mathcal W$. Then there exists a unique map $\mapp\phi{X_V}M$ such
that $\phi_W=\phi\circ R^W_V$ whenever $W\in\mathcal W$, defined as
follows. Given an $\mathcal A\in X_V$; since it is a completely prime
filter, we have that $W\in\mathcal A$ for some $W\in\mathcal W$; then
we put $\phi(\mathcal A)=\phi_W(\mathcal A\downarrow W)$, where
$\mathcal A\downarrow W=\{\,U\in\mathcal A:U\subset W\,\}$.)
We denote this cosheaf by $\Fun{Fil_0}(\fbasetop)$.

\begin{proposition}
Suppose that $\fbase$ is a $T_0$-space. The following are equivalent:
\begin{itemize}
\item[\textup{(i)}]     $\fbase$ is sober.
\item[\textup{(ii)}]    Each faded cosheaf (of sets) over $\fbase$ is
       isomorphic to the cosheaf $\Fun{CS}(f)$ for some map
       $\mapp fX\fbaseset$.
\item[\textup{(iii)}]   The cosheaf\/ $\Fun{Fil_0}(\fbasetop)$ is
       isomorphic to the cosheaf $\Fun{CS}(f)$ for some map
       $\mapp fX\fbaseset$.
\end{itemize}
\end{proposition}

\begin{proof} (i)$\implies$(ii) follows from the Theorem,
(ii)$\implies$(iii) is trivial.

(iii)$\implies$(i). Let $X_V$ be the components and $R^W_V$ inclusion maps
of $\Fun{Fil_0}(\fbasetop)$, defined by~\eqref{eq:4} and~\eqref{eq:5}.
Consider an arbitrary $\mathcal A\in X_\fbaseset$ and let $\fpinbase$ be
the point such that $\mathcal A\overset f\mapsto b$ (through the isomorphism
stated by foreground). Then for $W\in\fbasetop$ we have that $b\in W$
iff $\mathcal A\in\Img R^W_\fbaseset$ (since $\Img R^W_\fbaseset$ may be
identified with $\tube fW$---cf.~\eqref{eq:2}); and from other side
we obviously have $\Img R^W_\fbaseset=\{\,\mathcal B\in X_\fbaseset:
W\in\mathcal B\,\}$. Hence $\mathcal A=\Neighbs b\fbaseset$ and there is no
other $b'\in\fbaseset$ having this property, since $\fbase$ is a
$T_0$\nobreakdash-space. Therefore $\fbase$ is sober (and besides,
$f$ is one-to-one).
\end{proof}

\end{document}